\documentclass{amsart}
\usepackage[utf8]{inputenc}
\usepackage[T1]{fontenc}
\usepackage[english]{babel}
\usepackage{amsmath}
\usepackage{amsfonts}
\usepackage{amssymb}
\usepackage{graphicx}

\usepackage{xcolor}

\usepackage{amsthm}
\newtheorem{theorem}{Theorem}

\theoremstyle{definition}
\newtheorem{corollary}{Corollary}

\theoremstyle{definition}
\newtheorem{proposition}{Proposition}

\theoremstyle{definition}
\newtheorem{definition}{Definition}

\theoremstyle{definition}

\DeclareMathOperator{\Id}{Id}

\DeclareMathOperator{\tr}{tr}

\title[The Embedding Problem in Algebras with Involution]{The Embedding Problem in Algebras with Involution}
\author{Jonatan Andres Gomez Parada}
\address{Department of Mathematics \\ IMECC, UNICAMP \\ 
Sérgio Buarque de Holanda, 651, 13083-859 Campinas, SP, Brazil}
\email{j211980@dac.unicamp.br}
\thanks{This study was financed in part by the Coordena\c c\~ao de Aperfei\c coamento de Pessoal de N\'{\i}vel Superior - Brasil (CAPES) -
Finance Code 001.}

\date{}
\begin{document}

\begin{abstract}
Let $K$ be an algebraically closed field of characteristic zero, and let $A$ and $B$ be two simple algebras with involution over $K$. In this note we study the embedding problem for algebras with involution. More specifically, if  the algebra $A$ satisfies the polynomial identities with involution of the algebra $B$, we investigate whether there exists an embedding of $A$ into $B$ that preserves the involutions. 
\end{abstract}

\keywords{Polynomial identities;  identities with involution; matrix algebra}

\subjclass{16R10, 16W10, 16R50}

\maketitle

\section*{Introduction}
Let  $K$ be a field and let $K\langle X \rangle$ be the free associative algebra freely generated by the countable set of indeterminates $X$. We consider $K\langle X\rangle$ as the set of all polynomials in the non-commuting variables from the set $X$.  For a given associative algebra  $A$, a polynomial $f$ in $K\langle X\rangle$ is called a {\sl polynomial identity} of $A$ if $f$ evaluates to zero when its variables are substituted with arbitrary elements of $A$. An algebra that satisfies  a non-zero polynomial identity is called a {\sl PI-algebra}, and  the set of all identities for $A$, the T-ideal of $A$, is denoted by $\Id(A)$. Clearly $\Id(A)$ is an ideal, and moreover it is closed under all endomorphisms of $K\langle X \rangle$. The converse is also true: if an ideal in $K\langle X \rangle$ is closed under endomorphisms then it is the ideal of identities of some algebra $A$.

Let $A$ and $B$ be two PI-algebras over a field $K$. It is clear that if $A$ and $B$ are isomorphic then they satisfy the same set of polynomial identities. A natural question arises: If $A$ and $B$ satisfy the same set of polynomial identities, are $A$ and $B$  isomorphic? In general, the answer is no. For example, for any algebra $A$, it holds that $\Id(A) = \Id(A \oplus A)$. Also, considering $\mathbb{C}$ and $\mathbb{R}$, the algebras of complex numbers and  real numbers as $\mathbb{R}$-algebras, we have that  $\Id(\mathbb{C}) = \Id(\mathbb{R})$. Another, more sophisticated example arises when one considers the matrix algebra $M_2(\mathbb{R})$ and the real quaternion algebra $\mathbb{H}_{\mathbb{R}}$, which are not isomorphic because the latter is a division algebra while the former has zero divisors. On the other hand it is well known that  $\Id(\mathbb{H}_{\mathbb{R}}) = \Id(M_2(\mathbb{R}))$. In order to see the latter equality observe that if $A$ is a $K$-algebra and if $L$ is an extension field of $K$ then $A$ and $A_L=A\otimes_KL$ satisfy the same identities with coefficients in $K$ as long as $K$ is infinite. Then $M_2(\mathbb{R})\otimes_{\mathbb{R}}\mathbb{C} \cong \mathbb{H}_{\mathbb{R}}\otimes_{\mathbb{R}} \mathbb{C} \cong M_2(\mathbb{C})$.

The above examples show that algebras can satisfy the same set of polynomial identities but need not be isomorphic. In the first example, we have non-simple algebras; in the second, simple but non-central algebras; and in the third, central simple algebras over a non-algebraically closed field. 
Therefore, a natural restriction in our problem will be to focus  on central simple algebras over algebraically closed fields. Under this restriction, we obtain an immediate answer for finite-dimensional associative algebras. This answer is a direct consequence of the celebrated theorem of Amitsur and Levitzki.

Let $St_n$ be the {\sl standard polynomial of degree $n$} 		\[ St_n(x_1,\dots,x_n) = \sum_{\omega \in S_n}({\rm sign }\omega)x_{\omega(1)}\cdots x_{\omega(n)}, \] 		where $S_n$ is the symmetric group of degree $n$, and ${\rm sign} \omega$ is the sign of the permutation $\omega$. 

\textbf{Theorem (Amitsur–Levitzki).} 
\textit{
The matrix algebra $M_n(K)$ satisfies the standard identity of degree $2n$. It does not satisfy polynomial identities of lower degree and, up to a multiplicative constant, $St_{2n} = 0$ is the only multilinear polynomial identity of degree $2n$ for $M_n(K)$.}

Hence, if $A$ and $B$ are central simple algebras over the algebraically closed field $K$ then $A\cong M_n(K)$ and $B\cong M_m(K)$ for some $m$ and $n$. If $m<n$ then $A$ satisfies $St_{2m}$ while $B$ does not. Analogously if $n<m$ one gets an absurd, and hence $m=n$.

Other positive results, under the restrictions cited above, were given for:
\begin{itemize}
	\item Finite-dimensional Lie algebras by Kushkulei and Razmyslov (1983)
	\item Finite-dimensional Jordan algebras by Drensky and Racine (1992)
	\item Finite-dimensional algebras by Shestakov and Zaicev (2011).
\end{itemize}

The case of simple associative algebras graded by an abelian group was solved by Koshlukov and Zaicev (see \cite{koshlukov2010identities}), and the result was extended by Aljadeff and Haile to arbitrary groups (see \cite{aljadeff2014simple}). For finite-dimensional algebras graded by a semigroup, the positive answer was given by Bahturin and Yasumura (see \cite{bahturin2019graded}). In the latter paper very general results were obtained concerning the isomorphism of two algebraic systems provided they satisfy the same identical relations.

A more general problem is the {\it embedding problem}: Consider $A$ and $B$ two $K$-algebras, such that $A$ satisfies the polynomial identities of $B$ (and possibly some additional ones), we ask whether there exists an embedding of $A$ into $B$.

Initially, the embedding problem was considered in the following way:  Let $E$ be an affine PI-algebra, over a field $K$, and satisfying all identities of $M_n(K)$ for some $n$, can $A$ be embedded into $M_r(C)$ for some positive integer $r$ and some commutative $K$-algebra $C$? 

The first example of a PI-algebra satisfying all $n \times n$-matrix identities which is not embeddable in matrices was constructed by Amitsur in \cite{amitsur1970noncommutative}. Other examples about this problem were given in \cite{irving1986embeddability}. 

Now, concerning the question about the embedding problem for simple algebras over an algebraically closed field, there are several  positive results, for example, in the case of finite-dimensional associative algebras, and for algebras that are graded by an abelian group over a field of characteristic zero, \cite{david2012graded}.

Now let us turn our attention to algebras with an additional structure. The positive answer to the above problem in the case of algebras with trace was obtained by Procesi  in 1987, see \cite{procesi1987formal}, and for algebras with trace and involution by Berele in 1990, see \cite{berele1990matrices}. We recall that  an {\sl involution} on an algebra $A$ is an antiautomorphism of order two, that is, a linear map $* \colon A \to  A$ satisfying $(ab)^{*} = b^{*}a^{*}$ and $(a^{*})^{*} = a$, for every $a$, $b \in  A$. In the classification of involutions, these are called involutions of the first kind. A {\sl trace}  on an algebra $A$ is a linear map $\tr : A \to A$ such that for $a$, $b \in A$ one has $\tr(a)b = b\tr(a)$, $\tr(ab) = \tr(ba)$, and  $\tr(\tr(a)b) = \tr(a) \tr(b)$.

\textbf{Theorem (Procesi, 1987).} 
\textit{Let $S$ be a $K$-algebra with a trace function. If $S$ satisfies all of the trace identities of $M_n(K)$, then there is a trace preserving injection of $S$ into $M_n(A)$, for some commutative $K$-algebra $A$.}
		
\textbf{Theorem  (Berele, 1990).} 
\textit{Let $S$ be a $K$-algebra with  trace and involution. If $S$ satisfies all of the $*$-trace identities of $ (M_n(K),*) $, then there is a commutative $K$-algebra $A$ and a trace and involution preserving injection from $S$ to $(M_n(A),*)$.}

We draw the readers' attention to the following facts. Concerning the embedding problem, one may want to follow an argument based on the explicit form of the identities satisfied by the given algebras. But one stumbles in a significant and very hard problem: the identities satisfied even by central simple algebras are known in very few instances. The identities of $M_2(K)$ were described in \cite{razmyslov1973finite} in characteristic 0, and in \cite{koshlukov2001basis} in characteristic $p>2$. The identities of even $M_3(K)$ are rather far from our reach. The same applies for the concrete form of identities with involution in matrix algebras: these are known in the case of $M_2(K)$, see \cite{levchenko1982finite} e \cite{colombo2005identities}.

In this paper, we consider the embedding problem for simple algebras with involution over an algebraically closed field of characteristic different from $2$. We describe several cases where the answer to the embedding problem is positive.

\section{The Embedding Problem in Algebras with Involution}

Let $K$ be an algebraically closed field of characteristic different from $2$, and let $A$ be a $K$-algebra with involution $*$. We can decompose the algebra $A$ into the direct sum of the vector spaces of symmetric and skew-symmetric parts, in the following way. An element $a \in A$ is called {\sl symmetric} if $a^* = a$, and {\sl skew-symmetric} if $a^* = -a$. Thus, $a + a^*$ is symmetric and $a - a^*$ is skew-symmetric for any $a \in A$. Therefore, we have a decomposition $A = A^+ \oplus A^-$, where $A^+$ is the subspace of all symmetric elements, and $A^-$ is the subspace of all skew-symmetric elements of $A$. 

Now, we consider the polynomial identities of an algebra endowed with an involution. To this end, let $X = \{x_1, x_2, \ldots \}$ be a countable set of non-commutative variables, and denote $X^* = \{x_i, x_i^* \mid x_i \in X \}$. A $*$-action can be defined on monomials in $X^*$ by 
\[
(t_{i_1}\cdots t_{i_m})^* = t_{i_n}^* \cdots t_{i_1}^*, \quad t_{i} \in X^*, \, \text{ and } \, (x_j^*)^* = x_j \, \text{ for every } \, x_j \in X. 
\]
Let us consider $K\langle X, *\rangle = K\langle X^*  \rangle$, the free associative algebra freely generated by the set $X^*$. We call this algebra   the {\sl free associative algebra with involution}. Note that, in the case when the characteristic of $K$ is different from $2$, we can define $y_i = (x_i + x_i^*)/2$ and $z_i = (x_i - x_i^*)/2$ for each $i = 1$, 2,  \dots, and then consider $K\langle X, *\rangle = K\langle Y, Z \rangle = K\langle y_1, z_1, y_2, z_2, \ldots \rangle$ as generated by symmetric and skew-symmetric variables. The elements of $K\langle X, *\rangle$ will be called $*$-polynomials.

\begin{definition}
A $*$-polynomial $f\left( y_1,\dots,y_n,z_1,\dots,z_m \right) \in K\langle Y \cup Z \rangle$ is called a  {\sl $*$-polynomial identity} (or simply a {\sl $*$-identity}) of an algebra with involution $(A,*)$ whenever \[ f\left( u_1,\dots,u_n,v_1,\dots,v_m \right) = 0 \text{ for every } u_i \in A^+ \text{ and } v_j \in A^- . \] 
\end{definition}
Given an algebra with involution $(A,*)$ we denote by $\Id^*(A)$ or $T^*(A)$ the set of $*$-polynomial identities of  $A$. Clearly $T^*(A)$ is an ideal, closed under endomorphisms compatible with the involution.

 Recall that in the case of $M_n(K)$, the algebra of $n\times n$ matrices over $K$,  with $K$ an algebraically closed field of characteristic different from $2$, there are two types of involutions:
\begin{itemize}
	\item The {\it transpose involution}, denoted by $t$:
\[ \left(a_{ij}\right)^t = \left(a_{ji}\right), \; \text{ where } \; \left(a_{ij}\right) \in  M_n(K).\]
 	\item In the case when $n = 2k$ is even, the {\it symplectic involution} denoted by $s$ is defined as follows: \[ a^s = Ta^tT^{-1}, \text{ for all } a\in M_{2k}(K), \] where $T = \sum_{i=1}^{k}\left(e_{i,i+k} - e_{i+k,i}\right)$.
  
	That is, if $n = 2k$ and $B \in M_{2k}(K)$, we consider $B$ as a block matrix of size $k \times k$, and thus \[ \begin{bmatrix}
		R & S \\ P & Q \end{bmatrix}^s = \begin{bmatrix}
		Q^t & -S^t \\ - P^t & R^t
	\end{bmatrix}  \]
\end{itemize}

The classification of the finite dimensional simple algebras with involution, given below, is well known.

\begin{theorem}
	Let $K$ be an algebraically closed field. Every finite dimensional $*$-simple  $K$-algebra with involution $A$ is isomorphic, as a $*$-algebra, to one  of the following types:
	
	\begin{itemize}
		\item  $(M_k (K), t)$, the full matrix algebra with the transpose involution,
		\item $(M_k (K), s)$, the full matrix algebra with the symplectic involution, ($k\in 2\mathbb{Z}$),
		\item $(M_k (K) \oplus M_k (K)^{op}, ex)$, the direct product of the full matrix algebra and its opposite algebra with the exchange involution $*$.   
	\end{itemize}
\end{theorem}
Denote by $\left( M_n(K),*\right)$ the full matrix algebra with involution $*$.

A possible approach to the embedding problem can be based on considering the standard polynomials, as in the case of simple algebras. Hence we look for  the minimum degree of a standard polynomial that becomes an identity with involution for the matrix algebra.

By the Theorem of Amitsur and Levitzki, if $d\geq 2n$, the standard polynomial  $St_{d}(x_1, \dots , x_{2n-2})$ is a $*$-identity for $(M_{n}(K), *)$, where the  $x_i$'s can be symmetric or skew-symmetric variables, and $*$ any involution on $M_{n}(K)$. On the other hand, if $* = t$ (the transpose involution), the following theorem was obtained by Rowen (see \cite{rowen1980polynomial}).
 
 \begin{theorem} \label{th:skewvariables-t}
  Let $(M_{n}(K), t)$  be the algebra of $n \times n$ matrices endowed with the transpose involution. Then the polynomial $St_{2n-2}(z_1, \dots , z_{2n-2})$  is a standard $*$-identity of minimal degree in skew-symmetric variables.
 \end{theorem}
 
Let us consider now the algebras with involution $\left( M_n(K),t\right)$ and $\left( M_m(K),t\right)$ where $m > 1$, and suppose that ${\Id}\left( M_n(K),t\right) \subseteq {\Id}\left( M_m(K),t\right)$. By Theorem \ref{th:skewvariables-t}, we know that $St_{2n-2}(z_1, \dots , z_{2n-2})$, the standard polynomial of degree $2n-2$ in skew-symmetric variables is a $*$-identity for $\left( M_n(K),t\right)$.  Therefore, by assumption, we also have $St_{2n-2}(z_1, \dots , z_{2n-2}) \in \Id\left( M_m(K),t\right)$, and once again, by Theorem \ref{th:skewvariables-t}, we have $2n-2 \geq 2m-2$. Therefore, $n\geq m$.

 Thus, we have the injective $K$-homomorphism $\iota\colon M_m(K) \to M_n(K)$, given by \begin{equation}\label{eq:iota}
 \iota\left((a_{kl})\right) = (b_{kl}), \text{ where } b_{kl} = \begin{cases}
 a_{kl} & \text{ if } 1\leq k,l \leq m \\
 0 & \text{ otherwise.}
 \end{cases}  
 \end{equation}
 That is, $M_m(K)$ embeds into the upper left corner of $M_n(K)$. 
Therefore, $\iota\left((a_{kl})^t\right) = \iota\left((a_{kl})\right)^t$, and thus $\left( M_m(K),t\right)$ embeds into $\left( M_n(K),t\right)$, and we write $\left( M_m(K),t\right) \hookrightarrow  \left( M_n(K),t\right)$. 
Furthermore in the case when $m=1$, it is evident that $M_m(K)$ embeds into $M_n(K)$. 

Now, let us consider the case of the symplectic involution.

\begin{theorem}[\cite{StandardPolynomialMatrices}, Lemma 4.1] \label{skewvariables-s}
Let $(M_{2k}(K), s)$  be the algebra of $2k \times 2k$ matrices endowed with the symplectic involution. Then the polynomial $St_{4k}(x_1, \dots , x_{4k})$  is a standard $*$-identity of minimal degree in skew variables.
\end{theorem}

Let us now consider $\left( M_{2n}(K),s\right)$ and $\left( M_{2m}(K),s\right)$ with $m\geq 1$, and suppose that ${\rm Id}\left( M_{2n}(K),s\right) \subseteq {\rm Id}\left( M_{2m}(K),s\right)$. By Theorem \ref{skewvariables-s}, we have that $m\leq n$. Furthermore, if $m\leq n$, for $A \in M_m(K)$, we consider $\overline{A} = \iota(A) \in M_n(K)$ as in (\ref{eq:iota}), and we define \[ \varphi: M_{2m}(K) \to M_{2n}(K), \quad \begin{pmatrix} A & B \\ C & D \end{pmatrix} \mapsto \begin{pmatrix} \overline{A}  & \overline{B}  \\ \overline{C}  & \overline{D}  \end{pmatrix}. \] 

+ 
For example, \[ \varphi: M_2(K) \to M_4(K), \quad  \begin{pmatrix}
	a & b \\ c  & d \end{pmatrix} \mapsto  \begin{pmatrix}
	a & 0 & b & 0 \\ 0 & 0 & 0 & 0 \\ c & 0 & d & 0 \\ 0 & 0 & 0 & 0 \end{pmatrix}. \]
Then, $\varphi$ is an injective homomorphism that preserves the symplectic involution. Indeed, \[ \varphi \left( \begin{pmatrix} A & B \\ C & D \end{pmatrix}^s\right)  =  \begin{pmatrix} \overline{D^t} & \overline{-B^t} \\ \overline{-C^t} & \overline{A^t} \end{pmatrix}  = \begin{pmatrix} \overline{D}^t & - \overline{B}^t \\ -\overline{C}^t & \overline{A}^t \end{pmatrix}  =   \varphi \left( \begin{pmatrix} A & B \\ C & D \end{pmatrix}\right)^s . \]

Therefore, if ${\rm Id}\left( M_{2n}(K),s\right) \subseteq {\rm Id}\left( M_{2m}(K),s\right)$, then $\left( M_{2m}(K),s\right) \hookrightarrow \left( M_{2n}(K),s\right)$.

In the case when ${\rm Id}\left(  \left( M_{m}(K),t\right) \right)  \subseteq {\rm Id}\left(  \left( M_{2n}(K),s\right) \right)$, by Theorem \ref{th:skewvariables-t} and Theorem \ref{skewvariables-s}, we have $2m-2 \geq 4n$. Therefore, it follows that $m>m-1\geq 2n$, and as a consequence, $M_{2n}(K) \hookrightarrow M_{m}(K)$.

For standard identities in symmetric variables, we recall the following theorems.

\begin{theorem}[\cite{slin1979special}, Proposition 2] \label{symmetricvariables-t}
    Let $(M_{n}(K), t)$  be the algebra of $n \times n$ matrices endowed with the transpose involution. Then the polynomial $St_{2n}(x_1, \dots , x_{2n})$  is a standard $*$-identity of minimal degree in symmetric variables.
\end{theorem}

\begin{theorem}[\cite{rowen1982simple}, Theorem 3] \label{symmetricvariables-s}
	The standard polynomial $S_{4k-2}(y_1,\dots,y_{4k-2})$  is a $*$-identity of $(M_{2k}(K), s)$ in symmetric variables for every $k\geq 1$. 
\end{theorem}

Now, we consider some relations between $ \Id\left( \left( M_{2n}(K),s\right) \right)$ and $\Id\left( \left( M_{m}(K),t\right) \right)$. 

If $ \Id\left( \left( M_{2n}(K),s\right) \right) \subseteq \Id\left( \left( M_{m}(K),t\right) \right)$, by Theorem \ref{symmetricvariables-t} and Theorem \ref{symmetricvariables-s}, we have $4n - 2 \geq 2m$. Thus, $2n \geq m$ and $M_{m}(K) \hookrightarrow M_{2n}(K)$.\\

\begin{proposition}[\cite{becher2018involutions}, Proposition 4.4] \label{pro1}
Let $m$ be a positive integer. Then
\[
\Psi \colon (M_m(K), t) \to (M_{2m}(K), s), \qquad \alpha \mapsto  \begin{pmatrix}
		\alpha & \mathbf{0} \\ \mathbf{0} & \alpha \end{pmatrix}  
  \]
  is a homomorphism of $K$-algebras with involution.
\end{proposition}

\begin{corollary}
	${\rm Id}(M_{2m}(K),s) \subset  {\rm Id}(M_{m}(K), t) $.
\end{corollary}

\begin{proposition}[\cite{rowen1980polynomial}, Corollary 2.5.12] \label{pro2} 
	$(M_n(K),t)  \nsubseteq  (M_{2m}(K),s)$ for any  $m<n$.
\end{proposition}
Even though Proposition \ref{pro2}  specifically references algebras, in its proof, it was shown that if $m < n$, then ${\rm Id}(M_{2m}(K),s) \nsubseteq  {\rm Id}(M_{n}(K), t)$. 

From Proposition \ref{pro2}, if ${\rm Id}\left(  \left( M_{2m}(K),s \right) \right) \subseteq    {\rm Id}\left(  \left( M_{n}(K), t \right)\right)$, then $n \leq m$. Proposition \ref{pro1} implies that \[ \left( M_{n}(K),t\right) \hookrightarrow \left( M_{m}(K),t\right) \hookrightarrow \left( M_{2m}(K),s\right).\]
Based on the facts described above, we have the following result.
\begin{proposition}
Let us consider the algebras with involution $\left( M_{n_1}(K),*_1\right)$ and $\left( M_{n_2}(K),*_2\right)$, where $n_i$ is even if $*_i = s$.
	 \begin{itemize}
	 	\item[(a)] If  ${\rm Id}\left(  \left( M_{n_1}(K),*_1\right) \right) \subseteq    {\rm Id}\left(  \left( M_{n_2}(K),*_2\right)\right)$, then  $ M_{n_2}(K) \hookrightarrow M_{n_1}(K)$.
	 	
	 	\item[(b)] If  ${\rm Id}\left(  \left( M_{n_1}(K),*\right) \right) \subseteq    {\rm Id}\left(  \left( M_{n_2}(K),*\right) \right)$, then  $ \left( M_{n_2}(K),*\right)\hookrightarrow \left( M_{n_1}(K),*\right)$.

	  \item[(c)] If ${\rm Id}\left(  \left( M_{n_1}(K),s\right) \right) \subseteq    {\rm Id}\left(  \left( M_{n_2}(K),*\right) \right)$, then  $ \left( M_{n_2}(K),*\right)\hookrightarrow \left( M_{n_1}(K),s\right)$.
	 \end{itemize}
	\end{proposition}
As a consequence of the previous result, we have the following theorems:

\begin{theorem}
	Let $A$ and $B$ be two finite-dimensional central simple algebras with involution over the algebraically closed field $K$ of characteristic different from $2$,  $A$ with involution of orthogonal type and $A$ satisfying the identities with involution of the algebra $B$. Then, there exists an involution preserving embedding of $A$ into $B$.
\end{theorem}

\begin{theorem}
	Let $A$ and $B$ be two finite-dimensional central simple algebras with involution over the algebraically closed field $K$ of characteristic different from $2$,  $B$ with involution of symplectic type and $A$ satisfying the identities with involution of the algebra $B$. Then, there exists an embedding that preserves the involutions of $A$ into $B$.
\end{theorem}

Note that for $k\leq l$ one has that
\[  
(M_k (K) \times M_k (K)^{op}, ex) \to (M_l (K) \times M_l (K)^{op}, ex), \quad (A,B) \mapsto (\overline{A},\overline{B}) 
\]
and  
\[  (M_k (K) \times M_k (K)^{op}, ex) \to (M_{2k}(K),s), \quad (A,B) \mapsto  \begin{pmatrix}
	A & \mathbf{0} \\ \mathbf{0} & B^t \end{pmatrix} 
 \] 
are embeddings of $*$-algebras.

Therefore, if ${\rm Id}\left( (M_{2k}(K),s) \right) \subseteq \operatorname{Id}\left(  (M_l (K) \times M_l (K)^{op}, ex) \right)$ considering standard polynomials in symmetric variables, we have $4k - 2 \geq 2l$. Thus, $2k \geq l$, and we have an embedding that preserves the involutions from $(M_l (K) \times M_l (K)^{op}, ex)$ to $(M_{4k}(K),s)$.

\begin{corollary}
    Let $A$ be a finite-dimensional simple algebra with involution over the algebraically closed field $K$ of characteristic different from $2$ such that $A$ satisfies the identities with involution of the matrix algebra $(M_n(K), s)$.
    \begin{itemize}
        \item If $A$ is central, there exists an embedding of $A$ into $(M_{n}(K), s)$ that preserves the involutions. 
        \item If $A$ is not central, there exists an embedding of $A$ into $(M_{2n}(K), s)$ that preserves the involutions.
    \end{itemize}
\end{corollary}

\section*{Acknowledgements}
I would like to thank Plamen Koshlukov for suggesting the topic and for the helpful discussions.

\bibliographystyle{amsplain}
\bibliography{references}

\end{document}